\documentclass[11pt]{amsart}
\usepackage{amssymb,mathtools,fullpage,aliascnt,eucal,stmaryrd}
\usepackage[dvipsnames]{xcolor}
\usepackage[hypertexnames=false]{hyperref}
\usepackage[capitalize,nameinlink,noabbrev,nosort]{cleveref}
\hypersetup{
  colorlinks=true,
  allcolors=RedViolet
}
\DeclareFontEncoding{OT2}{}{}
\DeclareFontSubstitution{OT2}{cmr}{m}{n}
\DeclareFontFamily{OT2}{cmr}{\hyphenchar\font45 }
\DeclareFontShape{OT2}{cmr}{m}{n}{
<5><6><7><8><9>gen*wncyr
<10><10.95><12><14.4><17.28><20.74><24.88>wncyr10}{}
\DeclareFontShape{OT2}{cmr}{b}{n}{
<5><6><7><8><9>gen*wncyb
<10><10.95><12><14.4><17.28><20.74><24.88>wncyb10}{}
\DeclareMathAlphabet{\mathcyr}{OT2}{cmr}{m}{n}
\SetMathAlphabet{\mathcyr}{bold}{OT2}{cmr}{b}{n}
\newcounter{nameditem}

\newcommand{\nameditem}[1]{
  \item[(#1)]
  \renewcommand{\thenameditem}{(#1)}
  \refstepcounter{nameditem}
}
\theoremstyle{plain}
\newtheorem{theorem}{Theorem}[section]

\newaliascnt{proposition}{theorem}
\newtheorem{proposition}[proposition]{Proposition}
\aliascntresetthe{proposition}

\newaliascnt{conjecture}{theorem}
\newtheorem{conjecture}[conjecture]{Conjecture}
\aliascntresetthe{conjecture}
\crefname{conjecture}{Conjecture}{Conjectures}
\Crefname{conjecture}{Conjecture}{Conjectures}

\newaliascnt{lemma}{theorem}
\newtheorem{lemma}[lemma]{Lemma}
\aliascntresetthe{lemma}

\newaliascnt{corollary}{theorem}
\newtheorem{corollary}[corollary]{Corollary}
\aliascntresetthe{corollary}
\numberwithin{equation}{section}
\newcommand{\calA}{\mathcal{A}}
\newcommand{\rcf}{\mathrm{cf}}
\newcommand{\rCF}{\mathrm{CF}}
\newcommand{\dep}{\mathrm{dep}}
\newcommand{\rDe}{\mathrm{De}}
\newcommand{\be}{\boldsymbol{e}}
\DeclareMathOperator{\ev}{ev}
\newcommand{\calH}{\mathcal{H}}
\newcommand{\bh}{\boldsymbol{h}}
\newcommand{\bbI}{\mathbb{I}}
\newcommand{\calI}{\mathcal{I}}
\newcommand{\bk}{\boldsymbol{k}}
\newcommand{\rKZ}{\mathrm{KZ}}
\newcommand{\bl}{\boldsymbol{l}}
\newcommand{\bp}{\boldsymbol{p}}
\newcommand{\bbQ}{\mathbb{Q}}
\newcommand{\bbR}{\mathbb{R}}
\newcommand{\calR}{\mathcal{R}}
\newcommand{\reg}{\mathrm{reg}}
\newcommand{\calS}{\mathcal{S}}
\newcommand{\sh}{\mathbin{\mathcyr{sh}}}
\DeclareMathOperator{\Span}{span}
\newcommand{\bv}{\boldsymbol{v}}
\newcommand{\sfw}{\mathsf{w}}
\newcommand{\bbZ}{\mathbb{Z}}
\newcommand{\calZ}{\mathcal{Z}}
\newcommand{\rZ}{\mathrm{Z}}
\newcommand{\barZ}{\overline{\calZ}}
\newcommand{\fps}[1]{\langle\!\langle #1\rangle\!\rangle}
\newcommand{\llb}{\llbracket}
\newcommand{\rrb}{\rrbracket}
\begin{document}
\title{On the refined Kaneko--Zagier conjecture}
\author{Shin-ichiro Seki}
\address{Nagahama Institute of Bio-Science and Technology, 1266 Tamura-cho, Nagahama, Shiga, 526-0829, Japan}
\email{s\_seki@nagahama-i-bio.ac.jp}
\date{}
%--------------------------------------------------------------------------
\begin{abstract}
We extend our results on Zhao's conjectural basis to finite multiple zeta values modulo arbitrary prime powers and to the corresponding completed algebra.
As a consequence, the refined Kaneko--Zagier conjecture, which gives a precise description of relations modulo arbitrary prime powers in terms of ordinary multiple zeta values, follows from the dimension and direct-sum conjectures for ordinary and finite multiple zeta values.
We also prove Rosen's lifting conjecture under the dimension conjecture for finite multiple zeta values.
The proofs use confluence relations and $p$-adic associators and do not use motivic theory.
\end{abstract}
\maketitle
%--------------------------------------------------------------------------
\section{Introduction}
%--------------------------------------------------------------------------
\subsection{The Kaneko--Zagier conjecture}\label{subsec:KZconj}
%--------------------------------------------------------------------------
We follow the notation and conventions of \cite{Seki}.
An \emph{index} is a finite, possibly empty, tuple of positive integers.
For an index $\bk=(k_1,\dots,k_r)$, its \emph{weight} and \emph{depth} are defined by $\lvert\bk\rvert\coloneqq k_1+\cdots+k_r$ and $\dep(\bk)\coloneqq r$, respectively.
The empty index is denoted by $\varnothing$ and has weight and depth $0$.
An index is called \emph{admissible} if it is empty or its last entry is at least $2$.
We denote the set of all indices by $\bbI$ and set $\bbI_k=\{\bk\in\bbI:\lvert\bk\rvert=k\}$.

For an admissible index $\bk=(k_1,\dots, k_r)$, the \emph{multiple zeta value} $\zeta(\bk)$ is defined by
\[
\zeta(\bk)\coloneqq\sum_{0<n_1<\cdots<n_r}\frac{1}{n_1^{k_1}\cdots n_r^{k_r}}\in\bbR.
\]
We set $\zeta(\varnothing)\coloneqq1$ and use the same convention for the other zeta values below.
Let $\calZ_k$ denote the $\bbQ$-vector space generated by multiple zeta values of weight $k$, and put $\calZ\coloneqq\sum_{k=0}^{\infty}\calZ_k$ and $\calZ_k=0$ for $k<0$.
We write $\barZ\coloneqq\calZ/\zeta(2)\calZ$ and $\overline{\calZ_k}\coloneqq\calZ_k/\zeta(2)\calZ_{k-2}$.
For an index $\bk$, the \emph{finite multiple zeta value}
$\zeta^{}_{\calA}(\bk)$ is defined by
\[
\zeta^{}_{\calA}(\bk)\coloneqq
(\zeta^{}_{<p}(\bk)\bmod p)_p
\in\calA\coloneqq
\Biggl(\prod_{p:\text{ prime}}\bbZ/p\bbZ\Biggr)
\Bigg/\Biggl(\bigoplus_{p:\text{ prime}}\bbZ/p\bbZ\Biggr),
\]
where, for $\bk=(k_1,\dots,k_r)$, we set
\[
\zeta^{}_{<p}(\bk)\coloneqq
\sum_{0<n_1<\cdots<n_r<p}
\frac{1}{n_1^{k_1}\cdots n_r^{k_r}}.
\]
Set $\calZ_{\calA,k}\coloneqq\Span_{\bbQ}\{\zeta^{}_{\calA}(\bk):\bk\in\bbI_k\}$ and $\calZ_{\calA}\coloneqq\sum_{k=0}^{\infty}\calZ_{\calA,k}$.

For an arbitrary index $\bk=(k_1,\dots, k_r)$, let $\zeta^{\sh}(\bk)\in\calZ$ denote its shuffle-regularized multiple zeta value.
For each $0\leq i\leq r$, write $\bk_{[i]}=(k_1,\dots,k_i)$, $\bk^{[i]}=(k_{i+1},\dots,k_r)$, and $\overleftarrow{\bk}=(k_r,\dots,k_1)$.
The \emph{symmetric multiple zeta value} $\zeta^{}_{\calS}(\bk)$ is defined by
\[
\zeta^{}_{\calS}(\bk)\coloneqq\zeta^{\sh}_{\calS}(\bk)\bmod\zeta(2)\calZ,
\qquad
\zeta^{\sh}_{\calS}(\bk)\coloneqq\sum_{i=0}^{\dep(\bk)}(-1)^{\lvert\bk^{[i]}\rvert}\zeta^{\sh}(\bk_{[i]})\zeta^{\sh}(\overleftarrow{\bk^{[i]}})\in\calZ.
\]
Kaneko and Zagier \cite{KanekoZagier} conjectured that there exists an isomorphism of graded $\bbQ$-algebras
\[
\varphi_{\rKZ}\colon\calZ_{\calA}\overset{\sim}{\longrightarrow}\barZ
\]
such that $\varphi_{\rKZ}(\zeta^{}_{\calA}(\bk))=\zeta^{}_{\calS}(\bk)$.
The conjecture in \cite[Main Conjecture]{KanekoZagier} also includes a statement about the depth filtration, which we do not address in this paper.
See \cite[Remark~5.4]{Seki}.

Define a sequence of integers $(d_k)_{k\geq0}$ by
\[
\sum_{k=0}^{\infty}d_kX^k=\frac{1}{1-X^2-X^3},
\]
and set $d_{-1}=d_{-2}=0$ and $d_{-3}=1$.
For a positive integer $k$, let
\[
\bbI_k^{\rZ}\coloneqq\{(k_1,\dots,k_r,2,1)\in\bbI_k :r\geq0,\ k_1,\dots,k_r\in\{2,3\}\}.
\]
We also set $\bbI_0^{\rZ}\coloneqq\{\varnothing\}$.
The elements of these sets are called \emph{Zhao indices} in our convention.
We have $\#\bbI_k^{\rZ}=d_{k-3}$ for every $k\geq0$.

We consider the following four hypotheses concerning linear independence of zeta values, which will also be used in the main theorem of this paper:
\par\noindent
{\color{green!80!black}\rule{\linewidth}{0.6pt}}
\par
\textsc{Linear Independence Conjectures}
\par\vspace{2pt}
\begin{description}
\nameditem{Dim $\bbR$}\label{it:dim_MZV} $\dim_{\bbQ}\calZ_k\geq d_k$ for all $k\geq0$.
\nameditem{Dim $\calA$}\label{it:dim_FMZV} $\dim_{\bbQ}\calZ_{\calA,k}\geq d_{k-3}$ for all $k\geq3$.
\nameditem{DirSum $\bbR$}\label{it:direct_sum_MZV} $\displaystyle\calZ=\bigoplus_{k\geq0}\calZ_k$.
\nameditem{DirSum $\calA$}\label{it:direct_sum_FMZV} $\displaystyle\calZ_{\calA}=\bigoplus_{k\geq0}\calZ_{\calA,k}$.
\end{description}
\par\noindent
{\color{green!80!black}\rule{\linewidth}{0.6pt}}
\par\smallskip
Let $\calH\coloneqq\bbQ\langle e_0,e_1\rangle$ be the noncommutative polynomial ring in two variables $e_0$ and $e_1$ over $\bbQ$, and set $\calH^1\coloneqq\bbQ+e_1\calH$.
For an index $\bk=(k_1,\dots,k_r)$, define $\sfw(\bk)\coloneqq e_1e_0^{k_1-1}\cdots e_1e_0^{k_r-1}$.
Let $\calH^1_k$ be the subspace of $\calH^1$ generated by words of length $k$, with $\calH^1_0=\bbQ\cdot1$.
For each positive integer $k$, define the $\bbQ$-linear maps
$Z_{\calA,k}\colon\calH^1_k\longrightarrow\calZ_{\calA,k}$, $Z_{\calS,k}\colon\calH^1_k\longrightarrow\overline{\calZ_k}$ by sending $\sfw(\bk)$ to $\zeta^{}_{\calA}(\bk)$ and $\zeta^{}_{\calS}(\bk)$, respectively.
Note that, by abuse of notation, we write $\zeta^{}_{\calS}(\bk)\coloneqq\zeta^{\sh}_{\calS}(\bk)\bmod\zeta(2)\calZ_{k-2}$, with the intended meaning understood from the context.

The following theorem collects the results of \cite{Seki} relevant to the present paper.
\par\noindent
{\color{purple!80!black}\rule{\linewidth}{0.6pt}}
\par
\begin{theorem}[{\cite[Theorems~A, B, D, and E]{Seki}}]\label{thm:previous}
The following assertions hold.
\begin{enumerate}
\item\label{it:FMZV_basis} Every finite multiple zeta value is a $\bbQ$-linear combination of finite multiple zeta values with Zhao indices of the same weight.
\item\label{it:SMZV_basis} Every symmetric multiple zeta value is a $\bbQ$-linear combination of symmetric multiple zeta values with Zhao indices of the same weight.
\item\label{it:common_rel} For each positive integer $k$, there is an explicitly constructible subspace
\[
\calR_k\subset\ker Z_{\calA,k}\cap\ker Z_{\calS,k}, \qquad \dim_{\bbQ}\calR_k=2^{k-1}-d_{k-3}.
\]
A deterministic algorithm computes a basis of $\calR_k$ from $k$.
\item\label{it:KanekoZagier} Under \ref{it:dim_MZV}, \ref{it:dim_FMZV}, \ref{it:direct_sum_MZV}, and \ref{it:direct_sum_FMZV}, the Kaneko--Zagier conjecture holds.
\end{enumerate}
\end{theorem}
\par\vspace{-0.8\baselineskip}
\noindent
{\color{purple!80!black}\rule{\linewidth}{0.6pt}}
\par\smallskip
Analogues of \eqref{it:FMZV_basis} and \eqref{it:SMZV_basis}, with the Zhao indices replaced by a slightly different family of Hoffman-type indices, were proved independently and around the same time by Bachmann \cite{Bachmann}.
His proof applies Brown's method \cite{Brown2012} to a representation of motivic symmetric multiple zeta values obtained from Hirose's closed-path expression.
As explained in \cite{Seki}, both the existence of a common relation space of the stated dimension in \eqref{it:common_rel} and the assertion in \eqref{it:KanekoZagier} can be deduced from motivic results, notably \cite[Theorem~8.11]{Rosen}.
However, to the best of the author's knowledge, these statements have not been stated explicitly in the literature in this more elementary form, without reference to period maps, and may therefore be less familiar to non-specialists.
We also hope that some readers will find the proof in \cite{Seki} useful because it does not use the theory of mixed Tate motives.
In particular, the argument can be followed without going back to the underlying motivic work of Brown, Goncharov, and others, or to the construction of the category of mixed Tate motives over $\bbZ$.
We emphasize, however, that although our exposition is more elementary
and avoids the use of motives, the underlying arguments remain deeply
influenced by those developed by Brown \cite{Brown2012}.
Hirose's work \cite{Hirose} also provides the foundation that allows the arguments in \cite{Seki} to be developed without motivic theory.

The purpose of this note is to show that analogous results modulo arbitrary prime powers follow from \cite{Seki}, again without using motivic theory.
%--------------------------------------------------------------------------
\subsection{The refined Kaneko--Zagier conjecture}
%--------------------------------------------------------------------------
For a positive integer $n$, set
\[
\calA_n\coloneqq\Biggl(\prod_{p:\text{ prime}}\bbZ/p^n\bbZ\Biggr)\Bigg/\Biggl(\bigoplus_{p:\text{ prime}}\bbZ/p^n\bbZ\Biggr),
\qquad
\widehat{\calA}\coloneqq\varprojlim_n\calA_n.
\]
These are $\bbQ$-algebras, with rational numbers reduced at all primes not dividing their denominators.
Write $\bp\coloneqq(p\bmod p^n)_p\in\calA_n$ and also $\bp\coloneqq((p\bmod p^n)_p)_n\in\widehat{\calA}$.
For any index $\bk$, define
\[
\zeta^{}_{\calA_n}(\bk)\coloneqq(\zeta^{}_{<p}(\bk)\bmod p^n)_p\in\calA_n,\qquad\zeta_{\widehat{\calA}}(\bk)\coloneqq(\zeta^{}_{\calA_n}(\bk))_n\in\widehat{\calA}.
\]
These values were first introduced by Rosen \cite{Rosen2015}.
We use the topology on $\widehat{\calA}$ given by the kernels of its projections to $\calA_n$.
These kernels are $\bp^n\widehat{\calA}$.
Let $\calZ_{\calA_n}$ be the $\bbQ$-subalgebra of $\calA_n$ generated by $\bp$ and the values $\zeta^{}_{\calA_n}(\bk)$.
Following \cite[Section~4.1]{OnoSekiYamamoto2021}, put
\begin{equation}\label{eq:completed_algebra}
\calZ_{\widehat{\calA}}\coloneqq
\left\{\sum_{j=0}^{\infty}\bp^j
\sum_{\bk\in\bbI}a_{j,\bk}\zeta_{\widehat{\calA}}(\bk):
\begin{gathered}
a_{j,\bk}\in\bbQ,\\
\begin{aligned}
&\text{for each $j$, only finitely many}\\
&a_{j,\bk}\text{ are nonzero}
\end{aligned}
\end{gathered}\right\}.
\end{equation}
The series converge in $\widehat{\calA}$.
The stuffle product shows that this is a $\bbQ$-algebra.
Rosen proposed the following conjecture concerning relations among these zeta values.
\begin{conjecture}[Rosen {\cite[Conjecture~A]{Rosen2015}}]\label{conj:Rosen_lifting}
Let $n$ be a positive integer.
Then every homogeneous $\bbQ$-linear relation in $\calZ_{\calA_n}$ lifts to a $\bp$-adic relation of the same weight in $\calZ_{\widehat{\calA}}$.
Here $\bp$ is assigned weight $-1$.
\end{conjecture}
Takeyama and Tasaka observed that this conjecture follows from the injectivity of Rosen's $\calA$-valued period map; see \cite[Remark~3.2]{TakeyamaTasaka2023}.

To define the symmetric counterpart, let $\bk=(k_1,\dots,k_r)$ be an index and let $\bl=(l_1,\dots,l_r)\in\bbZ_{\geq0}^r$.
We write $\bk+\bl$ for their componentwise sum and set $\lvert\bl\rvert\coloneqq l_1+\cdots+l_r$.
Put
\[
b(\bk,\bl)\coloneqq\prod_{j=1}^{r}\binom{k_j+l_j-1}{l_j}.
\]
For the empty tuples, this product is $1$.
Let $t$ be an indeterminate.
For an index $\bk$, the \emph{$t$-adic symmetric multiple zeta value} $\zeta_{\widehat{\calS}}(\bk)$ is
\[
\zeta_{\widehat{\calS}}(\bk)\coloneqq\zeta^{\sh}_{\widehat{\calS}}(\bk)\bmod\zeta(2)\calZ\llb t\rrb\in\barZ\llb t\rrb,
\]
where
\[
\zeta^{\sh}_{\widehat{\calS}}(\bk)\coloneqq\sum_{i=0}^{\dep(\bk)}(-1)^{\lvert\bk^{[i]}\rvert}\zeta^{\sh}(\bk_{[i]})\sum_{\bl\in\bbZ_{\geq0}^{\dep(\bk)-i}}b(\bk^{[i]},\bl)\zeta^{\sh}(\overleftarrow{\bk^{[i]}+\bl})t^{\lvert\bl\rvert}\in\calZ\llb t\rrb.
\]
Its constant term is $\zeta^{}_{\calS}(\bk)$.
These values were first introduced by Jarossay \cite{JarossayAdjoint}.
Write $\zeta^{}_{\calS_n}(\bk)\coloneqq\zeta_{\widehat{\calS}}(\bk)\bmod t^n$.
The refinement of the Kaneko--Zagier conjecture has been considered in different formulations by Jarossay \cite{JarossayAdjoint} and Rosen \cite{Rosen2019}, among others.
Here we adopt the formulation of Ono--Seki--Yamamoto \cite[Conjecture~4.3]{OnoSekiYamamoto2021}.
We equip $\calZ_{\widehat{\calA}}$ with the subspace topology induced from $\widehat{\calA}$ and $\barZ\llb t\rrb$ with the $t$-adic topology.
\begin{conjecture}[Refined Kaneko--Zagier conjecture]\label{conj:refined_KZ}
There exists an isomorphism of topological $\bbQ$-algebras
\[
\widehat{\varphi}_{\rKZ}\colon\calZ_{\widehat{\calA}}
\overset{\sim}{\longrightarrow}\barZ\llb t\rrb
\]
such that $\widehat{\varphi}_{\rKZ}(\bp)=t$ and
$\widehat{\varphi}_{\rKZ}(\zeta_{\widehat{\calA}}(\bk))=\zeta_{\widehat{\calS}}(\bk)$ for every index $\bk$.
\end{conjecture}
For the assertions concerning weight, see \cref{thm:main} \ref{it:main_complete}.
Takeyama and Tasaka proved that the weightwise form of the Kaneko--Zagier conjecture (\cref{conj:refined_KZ}) implies the Kaneko--Zagier conjecture, provided that Rosen's lifting conjecture (\cref{conj:Rosen_lifting}) holds; see \cite[Proposition~3.6]{TakeyamaTasaka2023}.

For the finite-dimensional formulation, assign total weight $-1$ to both $\bp$ and $t$.
We use the conventions $\bbI_k=\varnothing$ and $\overline{\calZ_k}=0$ for $k<0$.
For an integer $k$, define
\begin{align}
\calZ_{\calA_n,k}&\coloneqq\Span_{\bbQ}\{\bp^j\zeta^{}_{\calA_n}(\bk):0\leq j<n,\ \bk\in\bbI_{k+j}\},\notag\\
\calZ_{\calS_n,k}&\coloneqq\bigoplus_{j=0}^{n-1}t^j \, \overline{\calZ_{k+j}}.\label{eq:symmetric_weight_space}
\end{align}
The second expression means that the coefficient of $t^j$ is taken in the indicated quotient.
There is a natural map from $\calZ_{\calS_n,k}$ to $\barZ[t]/(t^n)$.
Under \ref{it:direct_sum_MZV}, each $\calZ_{\calS_n,k}$ identifies with the homogeneous component of total weight $k$ in $\barZ[t]/(t^n)$.
When used in \eqref{eq:symmetric_weight_space}, the coefficients of $t^j\zeta^{}_{\calS_n}(\bk)$ are understood in these weightwise quotients.

The corresponding dimension conjecture is as follows (for $n=1$, this is Zhao's basis conjecture \cite[Conjecture~8.1 (ii)]{Zhao}).
\begin{conjecture}
For every positive integer $n$ and nonnegative integer $k$, we have
\[
\dim_{\bbQ}\calZ_{\calA_n,k}=\sum_{j=0}^{n-1}d_{k+j-3}.
\]
More precisely, the elements of $\{\bp^j\zeta^{}_{\calA_n}(\bl):0\leq j<n,\ \bl\in\bbI^{\rZ}_{k+j}\}$ form a basis of $\calZ_{\calA_n,k}$.
\end{conjecture}
A related conjecture was proposed by Zhao \cite[Conjecture~9.7]{Zhao}.
His space $\mathsf{FMZ}_{k,n}$ is generated by finite multiple zeta values whose indices have weight $k$, whereas $\calZ_{\calA_n,k}$ also includes the elements $\bp^j\zeta^{}_{\calA_n}(\bk)$ with $1\leq j<n$ and $\lvert\bk\rvert=k+j$.
For the relationship between these spaces, see \cref{sec:comparison}.

Set $\bbI^{\rZ}\coloneqq\bigcup_{k=0}^{\infty}\bbI^{\rZ}_k$.
For a positive integer $n$ and a nonnegative integer $k$, put
\[
\calH^{(n)}_k\coloneqq\bigoplus_{j=0}^{n-1}t^j\calH^1_{k+j}.
\]
Define $\bbQ$-linear maps $Z_{\calA_n,k}\colon\calH^{(n)}_k\longrightarrow\calZ_{\calA_n,k}$ and $Z_{\calS_n,k}\colon\calH^{(n)}_k\longrightarrow\calZ_{\calS_n,k}$ by $t^j\sfw(\bk)\mapsto\bp^j\zeta^{}_{\calA_n}(\bk)$ and
$t^j\sfw(\bk)\mapsto t^j\zeta^{}_{\calS_n}(\bk)$, respectively.
Our main result is as follows.
\newpage
\par\noindent
{\color{purple!80!black}\rule{\linewidth}{0.6pt}}
\par
\begin{theorem}\label{thm:main}
Let $n$ and $k$ be positive integers.
\begin{enumerate}
\nameditem{i}\label{it:main_gen_FMZV}
The space $\calZ_{\calA_n,k}$ is generated by $\{\bp^j\zeta^{}_{\calA_n}(\bl):0\leq j<n,\ \bl\in\bbI^{\rZ}_{k+j}\}$.
Moreover, every element of $\calZ_{\widehat{\calA}}$ admits an expansion
\[
\sum_{j=0}^{\infty}\bp^j\sum_{\bl\in\bbI^{\rZ}}a_{j,\bl}\zeta_{\widehat{\calA}}(\bl),
\]
where $a_{j,\bl}\in\bbQ$ and, for each $j$, only finitely many coefficients are nonzero.
\nameditem{ii}\label{it:main_gen_SMZV}
The space $\calZ_{\calS_n,k}$ is generated by $\{t^j\zeta^{}_{\calS_n}(\bl):0\leq j<n,\ \bl\in\bbI^{\rZ}_{k+j}\}$.
Moreover, every element of $\overline{\calZ}\llb t\rrb$ admits an expansion
\[
\sum_{j=0}^{\infty}t^j\sum_{\bl\in\bbI^{\rZ}}a_{j,\bl}\zeta_{\widehat{\calS}}(\bl),
\]
where $a_{j,\bl}\in\bbQ$ and, for each $j$, only finitely many coefficients are nonzero.
\nameditem{iii}\label{it:main_relations}
There is an explicitly constructible subspace
\[
\calR_k^{(n)}\subset\ker Z_{\calA_n,k}\cap\ker Z_{\calS_n,k},
\qquad
\dim_{\bbQ}\calR_k^{(n)}=2^{k-1}(2^n-1)-\sum_{j=0}^{n-1}d_{k+j-3},
\]
a basis of which can be computed by a deterministic algorithm from $k$ and $n$.
We have $\calR_k^{(1)}=\calR_k$ with $\calR_k$ as in \cref{thm:previous}, and reduction modulo $t^n$ maps $\calR_k^{(n+1)}$ onto $\calR_k^{(n)}$.
The construction of $\calR_k^{(n)}$ is given in \eqref{eq:def_relations_n} below.
\nameditem{iv}\label{it:main_complete}
Assume all four hypotheses \ref{it:dim_MZV}, \ref{it:dim_FMZV}, \ref{it:direct_sum_MZV}, and \ref{it:direct_sum_FMZV}.
Then the decompositions
\[
\calZ_{\calA_n}=\bigoplus_{k\in\bbZ}\calZ_{\calA_n,k},\qquad\barZ[t]/(t^n)=\bigoplus_{k\in\bbZ}\calZ_{\calS_n,k}
\]
give $\bbZ$-gradings, and there is an isomorphism of graded
$\bbQ$-algebras
\[
\varphi_{\rKZ}^{(n)}\colon\calZ_{\calA_n}\overset{\sim}{\longrightarrow}\barZ[t]/(t^n)
\]
such that $\varphi_{\rKZ}^{(n)}(\bp)=t$ and $\varphi_{\rKZ}^{(n)}(\zeta^{}_{\calA_n}(\bk))=\zeta^{}_{\calS_n}(\bk)$.
These isomorphisms are compatible with the natural reduction maps, and their inverse limit gives the isomorphism $\widehat{\varphi}_{\rKZ}$, that is, the refined Kaneko--Zagier conjecture (\cref{conj:refined_KZ}) holds.
\nameditem{v}\label{it:rosen_lifting}
Assume \ref{it:dim_FMZV}.
Then Rosen's lifting conjecture (\cref{conj:Rosen_lifting}) holds.
Moreover, there is a deterministic algorithm that, given a homogeneous $\bbQ$-linear relation in $\calZ_{\calA_n}$, computes any prescribed finite truncation of a $\bp$-adic lift of the same weight in $\calZ_{\widehat{\calA}}$.
\end{enumerate}
\end{theorem}
\par\vspace{-0.8\baselineskip}
\noindent
{\color{purple!80!black}\rule{\linewidth}{0.6pt}}
\par\smallskip
Except for the explicit dimension formula for $\calR_k^{(n)}$,
all assertions also hold for $k=0$, with $\calR_0\coloneqq\{0\}$.
In this case, we have
\[
\dim_{\bbQ}\calR_0^{(n)}=2^{n-1}-\sum_{j=0}^{n-1}d_{j-3}.
\]
The following dimension bound is a consequence of \ref{it:main_gen_FMZV}.
\begin{corollary}
For every positive integer $n$ and nonnegative integer $k$, we have
\begin{equation}\label{eq:upper_bound}
\dim_{\bbQ}\calZ_{\calA_n,k}\leq\sum_{j=0}^{n-1}d_{k+j-3}.
\end{equation}
\end{corollary}
The main results provide more than an upper bound.
Indeed, \ref{it:main_relations} gives
\[
\dim_{\bbQ}(\calH^{(n)}_k/\calR_k^{(n)})=\sum_{j=0}^{n-1}d_{k+j-3}
\]
and $\calR_k^{(n)}\subset\ker Z_{\calA_n,k}$.
Thus, we have succeeded in achieving the bound \eqref{eq:upper_bound} using an explicit family of relations.

A tuple of integers, whose entries need not be positive, is called a \emph{general index}.
The refined Kaneko--Zagier conjecture can be extended to general indices, and Ono--Yamamoto \cite[Theorem~1.1]{OnoYamamoto2023} showed that this extension follows from \cref{conj:refined_KZ}.
Consequently, the refined Kaneko--Zagier conjecture for general indices also holds under \ref{it:dim_MZV}, \ref{it:dim_FMZV}, \ref{it:direct_sum_MZV}, and \ref{it:direct_sum_FMZV}.

Among our main results, \ref{it:main_gen_FMZV}, \ref{it:main_gen_SMZV}, and the construction of the explicit relation spaces $\calR_k^{(n)}$ that reduce the upper bound for $\dim_{\mathbb Q}\mathcal Z_{\mathcal A_n,k}$ to the conjectured value are new.
The remaining assertions, however, are likely to have been known to a number of specialists since the work of Akagi--Hirose--Yasuda \cite{AkagiHiroseYasuda} and Rosen \cite{Rosen}.
We nevertheless include them here for the reasons explained at the end of \cref{subsec:KZconj}.
%--------------------------------------------------------------------------
\section{Preliminaries}
%--------------------------------------------------------------------------
\subsection{Formal multiple zeta values with respect to the confluence relations}\label{subsec:confluence-MZV}
%--------------------------------------------------------------------------
Let $\calH^0\coloneqq\bbQ+e_1\calH e_0$, and let $\calI_{\rCF}$ denote the shuffle ideal of confluence relations in the convention of Hirose--Sato \cite{HiroseSato2019}.
We use the algebra and formal multiple zeta values defined in \cite{Seki}:
\[
\calZ^{\rcf}\coloneqq\calH^0/\calI_{\rCF},\qquad
\zeta^{\rcf}(\bk)\coloneqq(-1)^{\dep(\bk)}\reg_{\sh}(\sfw(\bk))\bmod\calI_{\rCF}.
\]
Here $\reg_{\sh}\colon\calH\to\calH^0$ denotes the regularization map obtained by taking the constant term under the standard isomorphism $(\calH,\sh)\simeq(\calH^0,\sh)[e_0,e_1]$.
The sign is part of the convention for the confluence relations.
In particular, evaluation at the KZ associator sends $\zeta^{\rcf}(\bk)$ to $\zeta^{\sh}(\bk)$.
Let $\calZ^{\rcf}_k$ denote the $\bbQ$-vector space generated by formal multiple zeta values of weight $k$.
The algebra $\calZ^{\rcf}$ is graded by weight: $\calZ^{\rcf}=\bigoplus_{k=0}^{\infty}\calZ^{\rcf}_k$.
We set $\calZ^{\rcf}_k=0$ for $k<0$.
For any index $\bk$, put
\[
\zeta^{\rcf}_{\calS}(\bk)\coloneqq\sum_{i=0}^{\dep(\bk)}(-1)^{\lvert\bk^{[i]}\rvert}\zeta^{\rcf}(\bk_{[i]})\zeta^{\rcf}(\overleftarrow{\bk^{[i]}}).
\]
Write $\overline{\calZ^{\rcf}}\coloneqq\calZ^{\rcf}/\zeta^{\rcf}(2)\calZ^{\rcf}$, $\overline{\calZ^{\rcf}_k}\coloneqq\calZ^{\rcf}_k/\zeta^{\rcf}(2)\calZ^{\rcf}_{k-2}$.
Then $\overline{\calZ^{\rcf}}=\bigoplus_{k=0}^{\infty}\overline{\calZ^{\rcf}_k}$.

The following theorem restates \cite[Theorem~F]{Seki} in our notation, including the case of the empty index.
\begin{theorem}[{\cite[Theorem~F]{Seki}}]\label{thm:formal_basis}
For every nonnegative integer $k$, the elements of
\[
\{\zeta^{\rcf}(2)^a\zeta^{\rcf}_{\calS}(\bl):
a\in\bbZ \text{ with } 0\leq a\leq k/2, \ \bl\in\bbI^{\rZ}_{k-2a}\}
\]
form a basis of $\calZ^{\rcf}_k$.
Consequently, $\dim_{\bbQ}\calZ^{\rcf}_k=d_k$, and the elements of
\[
\{\zeta^{\rcf}_{\calS}(\bl)\bmod\zeta^{\rcf}(2) :\bl\in\bbI^{\rZ}_k\}
\]
form a basis of $\overline{\calZ^{\rcf}_k}$.
In particular, $\dim_{\bbQ}\overline{\calZ^{\rcf}_k}=d_{k-3}$.
\end{theorem}
As described in \cite[Section~1.2]{Seki}, we also obtain a deterministic algorithm for expansion in these bases.
First, shuffle regularization and the shuffle product express any prescribed homogeneous expression of weight $k$ in formal multiple zeta values as a rational linear combination of $\zeta^{\rcf}$-values for admissible indices.
Hirose's algorithm \cite[Section~1.3]{Hirose} then expresses this linear combination in the Hoffman basis.
Applying the same procedure to each element of the basis of $\calZ^{\rcf}_k$ in \cref{thm:formal_basis} gives an invertible rational transition matrix.
Using its inverse, we obtain the required expansion in that basis.
Discarding the terms with $a>0$ then gives the corresponding expansion in the basis of $\overline{\calZ^{\rcf}_k}$.
We will use this finite-dimensional computation repeatedly.

We denote the surjections induced by the real evaluation $\zeta^{\rcf}(\bk)\mapsto\zeta^{\sh}(\bk)$ by
\begin{equation}\label{eq:real_evaluations}
\ev_{\bbR,k}\colon\overline{\calZ^{\rcf}_k}\twoheadrightarrow\overline{\calZ_k},\qquad
\ev_{\bbR}\colon\overline{\calZ^{\rcf}}\twoheadrightarrow\barZ.
\end{equation}
The maps $\ev_{\bbR,k}$ are defined weight by weight, and $\ev_{\bbR}$ is an algebra homomorphism.
They send $\zeta^{\rcf}_{\calS}(\bk)\bmod\zeta^{\rcf}(2)$ to $\zeta^{}_{\calS}(\bk)$.
%--------------------------------------------------------------------------
\subsection{\texorpdfstring{Deligne's $p$-adic associator}{Deligne's p-adic associator}}\label{subsec:Deligne_ass}
%--------------------------------------------------------------------------
For a prime $p$, let $\Phi^p_{\rDe}\in\bbQ_p\fps{e_0,e_1}$ denote \emph{Deligne's $p$-adic associator}, with exactly the normalization used in \cite[Equation~(5.1)]{Seki}.
For every index $\bk$, we define \emph{Deligne's $p$-adic multiple zeta value} $\zeta_p^{\rDe}(\bk)$ by
\[
\zeta_p^{\rDe}(\bk)\coloneqq(-1)^{\dep(\bk)}[\overleftarrow{\sfw(\bk)}]\Phi^p_{\rDe}.
\]
Here $[w]\Phi$ denotes the coefficient of a word $w$ in $\Phi$, $\overleftarrow{w}$ is its reversal, and $\lvert w\rvert$ denotes its length.
We recall the $p$-adic inputs used in the proof of Theorem~A of \cite{Seki}.
Although we did not refer to this definition of $\zeta_p^{\rDe, \downarrow}(\bh)$ in \cite{Seki}, the expression given here is in fact the standard definition.
\begin{theorem}[Akagi--Hirose--Yasuda, Chatzistamatiou, Furusho, Jarossay, \"Unver]\label{thm:padic_inputs}
The following assertions hold.
\begin{enumerate}
\nameditem{i}\label{it:p-adic_associator}
The series $\Phi^p_{\rDe}$ is commutator group-like and satisfies the pentagon equation.
It defines an algebra homomorphism
\[
\ev_p\colon\overline{\calZ^{\rcf}}\longrightarrow\bbQ_p,
\qquad
\zeta^{\rcf}(\bk)\bmod\zeta^{\rcf}(2)\longmapsto\zeta_p^{\rDe}(\bk).
\]
\nameditem{ii}\label{it:padic_bound}
For $m\geq0$, let $(p\bbZ_p)^{[m]}$ be the ideal generated by $p^j/j!$ for all $j\geq m$.
For every word $w$,
\begin{equation}\label{eq:coefficient_bound}
[w]\Phi^p_{\rDe}\in p^{-\lvert w\rvert}(p\bbZ_p)^{[\lvert w\rvert]}.
\end{equation}
In particular, $\zeta_p^{\rDe}(\bk)\in\bbZ_p$ whenever $p>\lvert\bk\rvert$.
\nameditem{iii}\label{it:AHYconj}
For every index $\bk=(k_1,\dots,k_r)$, we have
\begin{equation}\label{eq:AHY_identity}
\zeta^{}_{<p}(\bk)=\sum_{i=0}^{\dep(\bk)}(-1)^{\lvert\bk^{[i]}\rvert}\zeta_p^{\rDe}(\bk_{[i]})\zeta_p^{\rDe,\downarrow}(\overleftarrow{\bk^{[i]}}),
\end{equation}
where, for every index $\bh$,
\begin{equation}\label{eq:down_arrow}
\zeta_p^{\rDe,\downarrow}(\bh)\coloneqq\sum_{j=0}^{\infty}p^j\sum_{\bl\in\bbZ_{\geq0}^{\dep(\bh)}, \ \lvert\bl\rvert=j}b(\bh,\bl)\zeta_p^{\rDe}(\bh+\bl).
\end{equation}
The series converges $p$-adically, and its value for the empty index is $1$.
\end{enumerate}
\end{theorem}
The assertion concerning $\Phi^p_{\rDe}$ in \ref{it:p-adic_associator} is due to \"Unver \cite{Unver2013}.
Combining this with $\zeta_p^{\rDe}(2)=0$ and applying the result of Furusho \cite{Furusho2022} shows that $\ev_p$ is well defined.
The bound \eqref{eq:coefficient_bound} follows from the integrality result of Chatzistamatiou \cite{Chatzistamatiou2017}.
Akagi--Hirose--Yasuda \cite{AkagiHiroseYasuda} have also announced an independent proof of this integrality result by a different method.
The identity in \ref{it:AHYconj} was conjectured by Akagi--Hirose--Yasuda and proved by Jarossay \cite{Jarossay}.
Any fixed element $z\in\overline{\calZ^{\rcf}}$ satisfies $\ev_p(z)\in\bbZ_p$ for all sufficiently large primes $p$.
\subsection{Formal $t$-adic symmetric multiple zeta values}
%--------------------------------------------------------------------------
For every index $\bk$, define $\zeta^{\rcf}_{\widehat{\calS}}(\bk)\in\overline{\calZ^{\rcf}}\llb t\rrb$ by
\begin{equation}\label{eq:formal_t_symmetric}
\zeta^{\rcf}_{\widehat{\calS}}(\bk)\coloneqq \sum_{i=0}^{\dep(\bk)}(-1)^{\lvert\bk^{[i]}\rvert}\sum_{\bl\in\bbZ_{\geq0}^{\dep(\bk)-i}}b(\bk^{[i]},\bl)\bigl(\zeta^{\rcf}(\bk_{[i]})\zeta^{\rcf}(\overleftarrow{\bk^{[i]}+\bl})\bmod\zeta^{\rcf}(2)\bigr)t^{\lvert\bl\rvert}.
\end{equation}
Its constant term is $\zeta^{\rcf}_{\calS}(\bk)\bmod{\zeta^{\rcf}(2)}$, and its coefficient of $t^j$ belongs to $\overline{\calZ^{\rcf}_{\lvert\bk\rvert+j}}$.
For each $j\geq 0$, $\ev_{\bbR}([t^j]\zeta^{\rcf}_{\widehat{\calS}}(\bk))=[t^j]\zeta_{\widehat{\calS}}(\bk)$.
\begin{proposition}\label{prop:modn-AkagiHiroseYasuda}
Let $\bk$ be an index and $n$ a positive integer.
If $p>\lvert\bk\rvert+n$, then
\begin{equation}\label{eq:modn-AkagiHiroseYasuda}
\zeta^{}_{<p}(\bk)\equiv
\sum_{j=0}^{n-1}p^j\ev_p\bigl([t^j]\zeta^{\rcf}_{\widehat{\calS}}(\bk)\bigr)
\pmod{p^n\bbZ_p}.
\end{equation}
Here, for each $0\leq j<n$, the value $\ev_p\bigl([t^j]\zeta^{\rcf}_{\widehat{\calS}}(\bk)\bigr)$ is a $p$-adic integer.
\end{proposition}
\begin{proof}
The assertion is immediate for the empty index.
For a nonempty index $\bh$ and $j\geq0$, \eqref{eq:coefficient_bound} and the integrality of $b(\bh,\bl)$ give
\[
p^j\sum_{\bl\in\bbZ_{\geq0}^{\dep(\bh)}, \ \lvert\bl\rvert=j}b(\bh,\bl)\zeta_p^{\rDe}(\bh+\bl)\in p^{-\lvert\bh\rvert}(p\bbZ_p)^{[\lvert\bh\rvert+j]}.
\]
If $p>\lvert\bh\rvert+n$ and $j\geq n$, then
\[
p^{-\lvert\bh\rvert}(p\bbZ_p)^{[\lvert\bh\rvert+j]}\subset p^{-\lvert\bh\rvert}(p\bbZ_p)^{[\lvert\bh\rvert+n]}=p^n\bbZ_p,
\]
where we have used $(p\bbZ_p)^{[m]}=p^m\bbZ_p$ for $p>m$.
Since the series in \eqref{eq:down_arrow} converges and $p^n\bbZ_p$ is closed, its tail with $j\geq n$ also belongs to $p^n\bbZ_p$.
Thus \eqref{eq:down_arrow} can be truncated at $j=n-1$ modulo $p^n\bbZ_p$.
The same truncation is immediate for the empty index.
Apply this truncation to each reversed suffix in \eqref{eq:AHY_identity}.
The factors $\zeta_p^{\rDe}(\bk_{[i]})$ are $p$-adic integers by \cref{thm:padic_inputs} \ref{it:padic_bound}, since $p>\lvert\bk\rvert$.
Replacing $\bl$ by $\overleftarrow{\bl}$ in each inner sum and using \eqref{eq:formal_t_symmetric}, we obtain \eqref{eq:modn-AkagiHiroseYasuda}.

For $0\leq j<n$, the coefficient
$[t^j]\zeta^{\rcf}_{\widehat{\calS}}(\bk)$
is, modulo $\zeta^{\rcf}(2)$, an integer linear combination of products of formal multiple zeta values, each product having total weight $\lvert\bk\rvert+j<p$.
Each factor therefore has $p$-adic integral evaluation by \cref{thm:padic_inputs} \ref{it:padic_bound}.
This proves the integrality assertion.
\end{proof}
\begin{proposition}\label{prop:evaluation_n}
For every positive integer $n$, there is an algebra homomorphism
\begin{equation}\label{eq:evaluation_n}
\psi_{\rKZ}^{\rcf, (n)}\colon\overline{\calZ^{\rcf}}[t]/(t^n)\longrightarrow\calA_n,
\qquad
\sum_{j=0}^{n-1}t^jz_j\longmapsto
\left(\sum_{j=0}^{n-1}p^j\ev_p(z_j)\bmod p^n\right)_p.
\end{equation}
It sends $t$ to $\bp$ and $\zeta^{\rcf}_{\widehat{\calS}}(\bk)\bmod t^n$ to $\zeta^{}_{\calA_n}(\bk)$.
The maps $\psi_{\rKZ}^{\rcf, (n)}$ are compatible with reduction and define a continuous algebra homomorphism
\[
\widehat{\psi}_{\rKZ}^{\rcf}\colon\overline{\calZ^{\rcf}}\llb t\rrb\longrightarrow\widehat{\calA}.
\]
\end{proposition}
\begin{proof}
For the finite collection of coefficients in \eqref{eq:evaluation_n}, all evaluations are integral outside a finite set of primes, by the final comment of \cref{subsec:Deligne_ass}.
The components at the excluded primes may be chosen arbitrarily.
Thus the right-hand side defines an element of $\calA_n$ independent of those choices.
Since each $\ev_p$ is an algebra homomorphism, so is $\psi_{\rKZ}^{\rcf, (n)}$.
The assertion about the symmetric series follows from \cref{prop:modn-AkagiHiroseYasuda}.
Compatibility of the maps is immediate from their definition, and the inverse limit defines $\widehat{\psi}^{\rcf}_{\rKZ}$.
The resulting map is continuous because its composite with the projection to $\calA_n$ factors through reduction modulo $t^n$ for every positive integer $n$.
\end{proof}
For a positive integer $n$ and a nonnegative integer $k$, write
\[
\overline{\calZ^{\rcf}_{\calS_n, k}}\coloneqq\bigoplus_{j=0}^{n-1}t^j\overline{\calZ^{\rcf}_{k+j}}\subset\overline{\calZ^{\rcf}}[t]/(t^n).
\]
Define the linear map $Z_{\calS_n, k}^{\rcf}\colon\calH_k^{(n)}\to\overline{\calZ^{\rcf}_{\calS_n, k}}$ by $t^j\sfw(\bk)\longmapsto t^j\zeta^{\rcf}_{\widehat{\calS}}(\bk)\bmod t^n$.
\begin{proposition}\label{prop:triangular_basis}
The elements of
\begin{equation}\label{eq:triangular_basis}
\{t^j\zeta^{\rcf}_{\widehat{\calS}}(\bl)\bmod t^n:
0\leq j<n,\ \bl\in\bbI^{\rZ}_{k+j}\}
\end{equation}
form a basis of $\overline{\calZ^{\rcf}_{\calS_n, k}}$.
In particular, $Z_{\calS_n, k}^{\rcf}$ is surjective and we have
\[
\dim_{\bbQ}\overline{\calZ^{\rcf}_{\calS_n, k}}=\sum_{j=0}^{n-1}d_{k+j-3}.
\]
\end{proposition}
\begin{proof}
The coefficient of $t^j$ in the element of \eqref{eq:triangular_basis} indexed by $(j,\bl)$ is $\zeta^{\rcf}_{\calS}(\bl)\bmod{\zeta^{\rcf}(2)}$, and its coefficients in lower powers of $t$ vanish.
By \cref{thm:formal_basis}, $\{\zeta^{\rcf}_{\calS}(\bl)\bmod{\zeta^{\rcf}(2)} : \bl\in\bbI^{\rZ}_{k+j}\}$ is a basis of $\overline{\calZ^{\rcf}_{k+j}}$.
Order both families by increasing $j$.
The coordinate matrix of the elements in \eqref{eq:triangular_basis} with respect to the basis $\{t^j(\zeta^{\rcf}_{\calS}(\bl)\bmod{\zeta^{\rcf}(2)}) : 0\leq j<n,\ \bl\in\bbI^{\rZ}_{k+j}\}$ is then block triangular with identity diagonal blocks.
It is invertible, proving the assertion.
\end{proof}
%--------------------------------------------------------------------------
\section{Proof of the main result}
%--------------------------------------------------------------------------
%--------------------------------------------------------------------------
\subsection{Construction of the space of relations $\calR_k^{(n)}$}\label{subsec:construction_of_rel}
%--------------------------------------------------------------------------
By uniqueness, expansions in the bases \eqref{eq:triangular_basis} are compatible with reduction.
Therefore, for every index $\bk$, there are unique rational coefficients $c_{\bk,\bl}^{(j)}$ such that
\begin{equation}\label{eq:infinite_formal_expansion}
\zeta^{\rcf}_{\widehat{\calS}}(\bk)=\sum_{j=0}^{\infty}t^j\sum_{\bl\in\bbI^{\rZ}_{\lvert\bk\rvert+j}}c_{\bk,\bl}^{(j)}\zeta^{\rcf}_{\widehat{\calS}}(\bl).
\end{equation}
We now describe a deterministic algorithm for computing these coefficients.
Suppose that the coefficients with superscript less than $j$ have been computed, and set
\[
z_j\coloneqq[t^j]\left(\zeta^{\rcf}_{\widehat{\calS}}(\bk)-\sum_{i=0}^{j-1}t^i\sum_{\bl\in\bbI^{\rZ}_{\lvert\bk\rvert+i}}c_{\bk,\bl}^{(i)}\zeta^{\rcf}_{\widehat{\calS}}(\bl)\right)\in\overline{\calZ^{\rcf}_{\lvert\bk\rvert+j}}.
\]
Expand $z_j$ in the basis $\{\zeta^{\rcf}_{\calS}(\bl)\bmod{\zeta^{\rcf}(2)} : \bl\in\bbI^{\rZ}_{\lvert\bk\rvert+j}\}$ using the algorithm recalled in \cref{subsec:confluence-MZV}.
The resulting coefficients are $c_{\bk,\bl}^{(j)}$.
Only finitely many coefficients of the series in \eqref{eq:formal_t_symmetric} are needed at each step.
If $\bk\in\bbI^{\rZ}$, uniqueness gives $c_{\bk,\bl}^{(0)}=\delta_{\bk,\bl}$ and $c_{\bk,\bl}^{(j)}=0$ for $j>0$, where $\delta_{\bk,\bl}$ is the Kronecker delta.

Let $n$ be a positive integer and $k$ a nonnegative integer.
For $0\leq j<n$ and $\bk\in\bbI_{k+j}\setminus\bbI^{\rZ}_{k+j}$, set
\begin{equation}\label{eq:def_relations_n}
\begin{split}
r_{j,\bk}^{(n)}&\coloneqq t^j\sfw(\bk)-\sum_{i=0}^{n-1-j}t^{j+i}
\sum_{\bl\in\bbI^{\rZ}_{k+j+i}}c_{\bk,\bl}^{(i)}\sfw(\bl),\\
\calR_k^{(n)}&\coloneqq\Span_{\bbQ}\{r_{j,\bk}^{(n)} : 0\leq j<n,\ \bk\in\bbI_{k+j}\setminus\bbI^{\rZ}_{k+j}\}.
\end{split}
\end{equation}
\begin{proposition}\label{prop:relations}
Let $n$ and $k$ be positive integers.
The displayed generators in \eqref{eq:def_relations_n} form a basis of $\ker Z_{\calS_n, k}^{\rcf}$.
Consequently, $\ker Z_{\calS_n, k}^{\rcf}=\calR_k^{(n)}$ and
\[
\dim_{\bbQ}\calR_k^{(n)}=2^{k-1}(2^n-1)-\sum_{j=0}^{n-1}d_{k+j-3}.
\]
\end{proposition}
\begin{proof}
The expansion \eqref{eq:infinite_formal_expansion} shows that each $r_{j,\bk}^{(n)}$ belongs to $\ker Z_{\calS_n,k}^{\rcf}$.
The coefficient of $t^j\sfw(\bk)$ in $r_{j,\bk}^{(n)}$ is $1$, whereas the coefficient of $t^{j'}\sfw(\bk')$ is $0$ whenever $0\leq j'<n$, $\bk'\in\bbI_{k+j'}\setminus\bbI^{\rZ}_{k+j'}$, and $(j,\bk)\neq(j',\bk')$.
Thus these elements are linearly independent.
By \cref{prop:triangular_basis}, their number equals the dimension of $\ker Z_{\calS_n,k}^{\rcf}$.
They therefore form a basis of that kernel. 
\end{proof}
\begin{lemma}\label{lem:two-composite-rep}
Let $n$ be a positive integer and $k$ a nonnegative integer.
Then
\[
Z_{\calA_n,k}=\psi_{\rKZ}^{\rcf,(n)}\circ Z_{\calS_n,k}^{\rcf}, \qquad
Z_{\calS_n,k}=\Biggl(\bigoplus_{j=0}^{n-1}t^j\ev_{\bbR,k+j}\Biggr)\circ Z_{\calS_n,k}^{\rcf}.
\]
\end{lemma}
\begin{proof}
This follows from the definitions and \cref{prop:evaluation_n}.
\end{proof}
\begin{proof}[Proof of \cref{thm:main} \ref{it:main_relations}]
By construction, $\calR_k^{(1)}=\calR_k$, and reduction modulo $t^n$ sends $r_{j,\bk}^{(n+1)}$ to $r_{j,\bk}^{(n)}$ if $j<n$ and to zero if $j=n$.
The remaining assertions follow from \cref{prop:relations,lem:two-composite-rep}.
\end{proof}
\begin{proof}[Proof of \cref{thm:main} \ref{it:main_gen_FMZV} and \ref{it:main_gen_SMZV}]
By \cref{prop:triangular_basis,lem:two-composite-rep}, we have $\psi_{\rKZ}^{\rcf,(n)}\bigl(\overline{\calZ^{\rcf}_{\calS_n,k}}\bigr)=\calZ_{\calA_n,k}$.
Thus the assertion concerning $\calZ_{\calA_n,k}$ in \ref{it:main_gen_FMZV} follows from \cref{prop:triangular_basis,prop:evaluation_n}.
Since each $\ev_{\bbR,k}$ is surjective, the assertion concerning $\calZ_{\calS_n,k}$ in \ref{it:main_gen_SMZV} follows from \cref{prop:triangular_basis}.

We next prove the assertions concerning $\calZ_{\widehat{\calA}}$ and $\barZ\llb t\rrb$.
First, let $f(t)\in\overline{\calZ^{\rcf}}\llb t\rrb$.
We construct rational coefficients $a_{j,\bl}$ recursively.
Suppose that the coefficients $a_{i,\bl}$ for $i<j$ have been chosen so that
\[
f(t)-\sum_{i=0}^{j-1}t^i\sum_{\bl\in\bbI^{\rZ}}a_{i,\bl}\zeta^{\rcf}_{\widehat{\calS}}(\bl)\in t^j\overline{\calZ^{\rcf}}\llb t\rrb,
\]
where each inner sum is finite.
Set
\[
z_j\coloneqq[t^j]\Biggl(f(t)-\sum_{i=0}^{j-1}t^i\sum_{\bl\in\bbI^{\rZ}}a_{i,\bl}\zeta^{\rcf}_{\widehat{\calS}}(\bl)\Biggr).
\]
Since $\overline{\calZ^{\rcf}}=\bigoplus_{k=0}^{\infty}\overline{\calZ^{\rcf}_k}$, the element $z_j$ has only finitely many nonzero weight components.
Applying \cref{thm:formal_basis} to each of these components, we obtain
\[
z_j=\sum_{\bl\in\bbI^{\rZ}}a_{j,\bl}\bigl(\zeta^{\rcf}_{\calS}(\bl)\bmod\zeta^{\rcf}(2)\bigr),
\]
with only finitely many nonzero coefficients.
The constant term of $\zeta^{\rcf}_{\widehat{\calS}}(\bl)$ is $\zeta^{\rcf}_{\calS}(\bl)\bmod\zeta^{\rcf}(2)$.
Thus subtracting $t^j\sum_{\bl\in\bbI^{\rZ}}a_{j,\bl}\zeta^{\rcf}_{\widehat{\calS}}(\bl)$ from the residual also makes its coefficient of $t^j$ vanish.
Starting with $j=0$ and repeating this procedure gives
\begin{equation}\label{eq:exp_f}
f(t)=\sum_{j=0}^{\infty}t^j\sum_{\bl\in\bbI^{\rZ}}a_{j,\bl}\zeta^{\rcf}_{\widehat{\calS}}(\bl),
\end{equation}
where, for each $j$, only finitely many coefficients are nonzero.
The equality holds because the two sides agree modulo $t^n$ for every positive integer $n$.
By \cref{prop:evaluation_n}, the continuous map $\widehat{\psi}_{\rKZ}^{\rcf}$ sends $t$ to $\bp$ and $\zeta^{\rcf}_{\widehat{\calS}}(\bk)$ to $\zeta_{\widehat{\calA}}(\bk)$.
Applying this map to the preceding expansion gives
\begin{equation}\label{eq:exp_psi_f}
\widehat{\psi}_{\rKZ}^{\rcf}(f(t))=\sum_{j=0}^{\infty}\bp^j\sum_{\bl\in\bbI^{\rZ}}a_{j,\bl}\zeta_{\widehat{\calA}}(\bl).
\end{equation}
Thus the image of $\widehat{\psi}_{\rKZ}^{\rcf}$ is contained in $\calZ_{\widehat{\calA}}$.
The reverse inclusion is immediate from \eqref{eq:completed_algebra}.
We have therefore proved
\begin{equation}\label{eq:image_of_E}
\widehat{\psi}_{\rKZ}^{\rcf}\bigl(\overline{\calZ^{\rcf}}\llb t\rrb\bigr)=\calZ_{\widehat{\calA}}.
\end{equation}
Together with the expansion \eqref{eq:exp_psi_f}, this proves the assertion concerning $\calZ_{\widehat{\calA}}$ in \ref{it:main_gen_FMZV}.
Finally, the surjection $\ev_{\bbR}\colon\overline{\calZ^{\rcf}}\twoheadrightarrow\barZ$ extends coefficientwise to a surjection $\ev_{\bbR}\colon\overline{\calZ^{\rcf}}\llb t\rrb\twoheadrightarrow\barZ\llb t\rrb$.
It sends $\zeta^{\rcf}_{\widehat{\calS}}(\bl)$ to $\zeta_{\widehat{\calS}}(\bl)$.
Applying it to the expansion \eqref{eq:exp_f} proves the assertion concerning $\barZ\llb t\rrb$ in \ref{it:main_gen_SMZV}.
\end{proof}
The inductive construction in this section is similar to those in Jarossay \cite[Proposition~5.5]{Jarossay2020} and Ono--Seki--Yamamoto \cite[Proposition~4.4]{OnoSekiYamamoto2021}.
The only difference is that we use the stronger result on expansion in the Zhao basis as input, in place of Yasuda's theorem \cite[Theorem~6.4]{Yasuda2016}.
%--------------------------------------------------------------------------
\subsection{Linear independence and the isomorphism between $\calZ_{\calA_n}$ and $\barZ[t]/(t^n)$}
%--------------------------------------------------------------------------
We first record the elementary cancellation property needed in the argument.
\begin{lemma}\label{lem:cancellation}
Let $a$ and $n$ be integers with $0\leq a<n$, and let $x\in\calA_n$.
If $\bp^ax=0$, then the image of $x$ in $\calA_{n-a}$ is zero.
In particular, if $\bp^ax=0$, then the image of $x$ in $\calA$ is zero.
\end{lemma}
\begin{proof}
Choose representatives $x_p\in\bbZ/p^n\bbZ$.
The equality $\bp^ax=0$ means that $p^ax_p\equiv0\pmod{p^n}$ for all sufficiently large primes $p$.
For those primes, $x_p\equiv0\pmod{p^{n-a}}$.
This is exactly the asserted vanishing in $\calA_{n-a}$.
\end{proof}
Let $\psi^{\rcf, (1)}_{\rKZ}\colon\overline{\calZ^{\rcf}}\to\calA$ be the evaluation in \cref{prop:evaluation_n} for $n=1$.
It sends $\zeta^{\rcf}_{\calS}(\bk)\bmod{\zeta^{\rcf}(2)}$ to $\zeta^{}_{\calA}(\bk)$.
\begin{lemma}\label{lem:residue_injectivity}
Assume \ref{it:dim_FMZV}.
Then $\psi^{\rcf, (1)}_{\rKZ}$ is injective on each $\overline{\calZ^{\rcf}_k}$.
If \ref{it:direct_sum_FMZV} also holds, then $\psi^{\rcf, (1)}_{\rKZ}$ is injective on $\overline{\calZ^{\rcf}}$.
\end{lemma}
\begin{proof}
\cref{thm:previous} \eqref{it:FMZV_basis} and \ref{it:dim_FMZV} imply that the values with indices in $\bbI^{\rZ}_k$ form a basis of $\calZ_{\calA,k}$ for every $k\geq3$.
For $k=0$, the assertion follows from $\psi^{\rcf, (1)}_{\rKZ}(1)=1$; for $k=1,2$, both spaces are zero.
Comparing with \cref{thm:formal_basis} proves injectivity on each weight.
Under \ref{it:direct_sum_FMZV}, the images of distinct weights have no nontrivial finite linear relation, proving injectivity on their direct sum.
\end{proof}

\begin{proposition}\label{prop:finite_injectivity}
Under \ref{it:dim_FMZV}, the restriction of $\psi^{\rcf, (n)}_{\rKZ}$ to $\overline{\calZ^{\rcf}_{\calS_n, k}}$ is injective for every nonnegative integer $k$ and every positive integer $n$.
Under both \ref{it:dim_FMZV} and \ref{it:direct_sum_FMZV}, $\psi^{\rcf, (n)}_{\rKZ}$ is injective on all of $\overline{\calZ^{\rcf}}[t]/(t^n)$.
\end{proposition}
\begin{proof}
Let $z=\sum_{j=0}^{n-1}t^jz_j$ be a nonzero element in $\overline{\calZ^{\rcf}_{\calS_n, k}}$, and let $a$ be the smallest index for which $z_a\neq0$.
Write $z=t^az'$ in $\overline{\calZ^{\rcf}}[t]/(t^n)$, where the constant term of $z'$ is $z_a$.
If $\psi^{\rcf, (n)}_{\rKZ}(z)=0$, then $\bp^a\psi^{\rcf, (n)}_{\rKZ}(z')=0$.
By \cref{lem:cancellation}, $\psi^{\rcf, (n)}_{\rKZ}(z')\bmod{\bp}$ is zero, so $\psi^{\rcf, (1)}_{\rKZ}(z_a)=0$.
In the first assertion, $z_a\in\overline{\calZ^{\rcf}_{k+a}}$, and the first part of \cref{lem:residue_injectivity} gives a contradiction.
For the second assertion, start with a nonzero element $z\in\overline{\calZ^{\rcf}}[t]/(t^n)$ and use the second part of \cref{lem:residue_injectivity} in the same argument.
\end{proof}
\begin{lemma}\label{lem:real_injectivity}
Under \ref{it:dim_MZV}, every $\ev_{\bbR,k}\colon\overline{\calZ^{\rcf}_k}\twoheadrightarrow\overline{\calZ_k}$ in \eqref{eq:real_evaluations} is an isomorphism.
If \ref{it:direct_sum_MZV} also holds, then $\ev_{\bbR}\colon\overline{\calZ^{\rcf}}\twoheadrightarrow\barZ$ is an isomorphism of graded $\bbQ$-algebras.
\end{lemma}
\begin{proof}
The formal upper bound in \cref{thm:formal_basis} and \ref{it:dim_MZV} give $\dim_{\bbQ}\calZ_k=d_k$.
Multiplication by the nonzero real number $\zeta(2)$ is injective, and hence $\dim_{\bbQ}\overline{\calZ_k}=d_k-d_{k-2}=d_{k-3}$.
Thus the surjection $\ev_{\bbR,k}$ has equal-dimensional source and target.

Under \ref{it:direct_sum_MZV}, the natural map $\bigoplus_{k=0}^{\infty}\overline{\calZ_k}\to\barZ$ is an isomorphism; see \cite[Lemma~5.3]{Seki}.
Combining this with the preceding weightwise isomorphisms proves the second assertion.
The resulting isomorphism preserves weight by construction.
\end{proof}

\begin{proof}[Proof of the first part of \cref{thm:main} \ref{it:main_complete}]
Assume \ref{it:dim_FMZV}.
By \cref{prop:finite_injectivity}, $\psi^{\rcf, (n)}_{\rKZ}$ identifies $\overline{\calZ^{\rcf}_{\calS_n, k}}$ with $\calZ_{\calA_n,k}$.
Under \ref{it:dim_MZV}, \cref{lem:real_injectivity} identifies the same formal space with $\calZ_{\calS_n,k}$.

Now assume the two direct-sum hypotheses as well.
Since $\psi_{\rKZ}^{\rcf,(n)}\bigl(\overline{\calZ^{\rcf}_{\calS_n,k}}\bigr)=\calZ_{\calA_n,k}$ for every $k\geq0$ and $\psi_{\rKZ}^{\rcf,(n)}(t)=\bp$, the image of $\psi_{\rKZ}^{\rcf,(n)}$ is $\calZ_{\calA_n}$.
The second assertion of \cref{prop:finite_injectivity} therefore gives an algebra isomorphism $\psi^{\rcf, (n)}_{\rKZ}\colon\overline{\calZ^{\rcf}}[t]/(t^n)\xrightarrow{\sim}\calZ_{\calA_n}$.
Assigning weight $-1$ to $t$, the source is a graded algebra
whose homogeneous component of total weight $k\in\bbZ$ is
\[
\bigoplus_{j=0}^{n-1}t^j\overline{\calZ^{\rcf}_{k+j}}.
\]
For $k\geq0$, its image under $\psi_{\rKZ}^{\rcf,(n)}$ is
$\calZ_{\calA_n,k}$, as already shown.
For $k=-h<0$, the same assertion follows from
\[
\bigoplus_{j=0}^{n-1}t^j\overline{\calZ^{\rcf}_{j-h}}
=t^h\overline{\calZ^{\rcf}_{\calS_n,0}},
\qquad
\calZ_{\calA_n,-h}=\bp^h\calZ_{\calA_n,0}.
\]
Since $\psi_{\rKZ}^{\rcf,(n)}$ is injective, these images give
a direct-sum decomposition of $\calZ_{\calA_n}$.
This decomposition defines a grading for which
$\psi_{\rKZ}^{\rcf,(n)}$ is a graded algebra isomorphism.

The graded algebra isomorphism $\ev_{\bbR}$ of \cref{lem:real_injectivity}, extended coefficientwise and fixing $t$, identifies its source with $\barZ[t]/(t^n)$ as a graded algebra.
Consequently,
\begin{equation}\label{eq:finite_isomorphism}
\varphi_{\rKZ}^{(n)}=(\ev_{\bbR}\bmod{t^n})\circ(\psi^{\rcf,(n)}_{\rKZ})^{-1}
\end{equation}
has all the required properties.
Its compatibility with reduction follows from that of the two evaluations.
\end{proof}
%--------------------------------------------------------------------------
\subsection{Passage to the completion}
%--------------------------------------------------------------------------
We complete the proof by checking the topology and the inverse limit.
This also explains why no further hypothesis about relations modulo higher powers of $p$ is necessary.
\begin{proposition}\label{prop:completed_evaluation}
Assume \ref{it:dim_FMZV} and \ref{it:direct_sum_FMZV}.
Then the map $\widehat{\psi}^{\rcf}_{\rKZ}$ in \cref{prop:evaluation_n} induces an isomorphism of topological $\bbQ$-algebras
\[
\widehat{\psi}^{\rcf}_{\rKZ}\colon\overline{\calZ^{\rcf}}\llb t\rrb\overset{\sim}{\longrightarrow}\calZ_{\widehat{\calA}}.
\]
Moreover, for every positive integer $n$, we have
\begin{equation}\label{eq:strictness}
\calZ_{\widehat{\calA}}\cap\bp^n\widehat{\calA}=\bp^n\calZ_{\widehat{\calA}},\qquad\calZ_{\widehat{\calA}}/\bp^n\calZ_{\widehat{\calA}}\simeq\calZ_{\calA_n},
\end{equation}
and $\calZ_{\widehat{\calA}}\simeq\varprojlim_n\calZ_{\calA_n}$.
\end{proposition}
\begin{proof}
Surjectivity follows from \eqref{eq:image_of_E}.
If $\widehat{\psi}^{\rcf}_{\rKZ}(f(t))=0$, then $\psi^{\rcf, (n)}_{\rKZ}(f(t)\bmod{t^n})=0$ for every $n$.
By \cref{prop:finite_injectivity}, $f(t)\bmod{t^n}=0$ for every $n$, so $f(t)=0$.
The same argument shows that $\widehat{\psi}^{\rcf}_{\rKZ}(f(t))\in\ker \, (\widehat{\calA}\to\calA_n)\Longleftrightarrow f(t)\in t^n\overline{\calZ^{\rcf}}\llb t\rrb$.
Thus $\widehat{\psi}^{\rcf}_{\rKZ}$ and its inverse are continuous.
Since $\widehat{\psi}^{\rcf}_{\rKZ}(t)=\bp$, this equivalence also proves the first equality in \eqref{eq:strictness}, and the second follows from the image of $\psi^{\rcf, (n)}_{\rKZ}$.
Finally, taking the inverse limit of the compatible isomorphisms $\psi^{\rcf,(n)}_{\rKZ}$ gives $\calZ_{\widehat{\calA}}\simeq\varprojlim_n\calZ_{\calA_n}$.
This proves the last assertion.
\end{proof}
\begin{proof}[Completion of the proof of \cref{thm:main}\,\ref{it:main_complete}]
Under all four hypotheses, extend the algebra isomorphism $\ev_{\bbR}\colon\overline{\calZ^{\rcf}}\to\barZ$ coefficientwise and compose it with the inverse of the map in \cref{prop:completed_evaluation}.
This gives a topological algebra isomorphism
\[
\widehat{\varphi}_{\rKZ}\coloneqq\ev_{\bbR}\circ(\widehat{\psi}^{\rcf}_{\rKZ})^{-1}\colon\calZ_{\widehat{\calA}}\overset{\sim}{\longrightarrow}\barZ\llb t\rrb.
\]
By \cref{prop:evaluation_n}, the map $\widehat{\psi}^{\rcf}_{\rKZ}$ sends $\zeta^{\rcf}_{\widehat{\calS}}(\bk)$ to $\zeta_{\widehat{\calA}}(\bk)$.
Real evaluation sends the same formal series to $\zeta_{\widehat{\calS}}(\bk)$.
Also, $\widehat{\psi}^{\rcf}_{\rKZ}(t)=\bp$ and real evaluation fixes $t$.
Hence $\widehat{\varphi}_{\rKZ}$ has precisely the properties in \cref{conj:refined_KZ}.
Its reduction at each truncation is \eqref{eq:finite_isomorphism}.
Under the natural identifications with the inverse limits,
we therefore have
\[
\widehat{\varphi}_{\rKZ}
=\varprojlim_n\varphi_{\rKZ}^{(n)}.
\]
Since each $\varphi_{\rKZ}^{(n)}$ is graded, the completed
isomorphism is compatible with total weight at each truncation.
\end{proof}
%--------------------------------------------------------------------------
\subsection{Rosen's lifting conjecture}
%--------------------------------------------------------------------------
We recall Rosen's formulation \cite{Rosen2015}, with the order of indices adapted to our convention.
Let $\widehat{\calH}\coloneqq\bbQ\fps{e_0,e_1}$ be the noncommutative formal power series ring in $e_0$ and $e_1$ over $\bbQ$, and set $\widehat{\calH}^1\coloneqq\bbQ+e_1\widehat{\calH}$.
Every element of $\widehat{\calH}^1$ can be written uniquely as
$\sum_{\bk\in\bbI}a_{\bk}\sfw(\bk)$ with $a_{\bk}\in\bbQ$.
For each positive integer $n$, set
\[
\calI_n\coloneqq\left\{\sum_{\bk\in\bbI}a_{\bk}\sfw(\bk)\in\widehat{\calH}^1 : a_{\bk}=0\text{ whenever }\lvert\bk\rvert<n\right\}.
\]
We equip $\widehat{\calH}^1$ with the topology for which the subspaces $\calI_n$ form a neighborhood basis of zero.

Define the continuous $\bbQ$-linear maps
\[
Z_{\widehat{\calA}}\colon\widehat{\calH}^1\longrightarrow\widehat{\calA},\qquad\sum_{\bk\in\bbI}a_{\bk}\sfw(\bk)\longmapsto\sum_{\bk\in\bbI}a_{\bk}\bp^{\lvert\bk\rvert}\zeta_{\widehat{\calA}}(\bk)
\]
and
\[
Z_{\calA_n}\colon\widehat{\calH}^1\longrightarrow\calA_n,\qquad\sum_{\bk\in\bbI}a_{\bk}\sfw(\bk)\longmapsto\sum_{\bk\in\bbI, \ \lvert\bk\rvert<n}a_{\bk}\bp^{\lvert\bk\rvert}\zeta^{}_{\calA_n}(\bk).
\]
The series defining $Z_{\widehat{\calA}}$ is summed in order of increasing weight and converges in $\widehat{\calA}$, since there are only finitely many indices of each weight and each summand belongs to $\bp^{\lvert\bk\rvert}\widehat{\calA}$.
The map $Z_{\calA_n}$ is the composite of $Z_{\widehat{\calA}}$ with the projection to $\calA_n$.
In particular, $\calI_n\subset\ker Z_{\calA_n}$.
Since $Z_{\widehat{\calA}}(\calI_n)\subset\bp^n\widehat{\calA}$ for every positive integer $n$, both maps are continuous, where $\calA_n$ is given the discrete topology.

In this notation, Rosen's lifting conjecture (\cref{conj:Rosen_lifting}) takes the form
\begin{equation}\label{eq:rosen_kernel_equality}
\ker Z_{\calA_n}=\ker Z_{\widehat{\calA}}+\calI_n
\end{equation}
for all $n\geq 1$.
Indeed, membership in $\ker Z_{\calA_n}$ means that
\[
\sum_{\bk\in\bbI, \ \lvert\bk\rvert<n}a_{\bk}p^{\lvert\bk\rvert}\zeta^{}_{<p}(\bk)\equiv0\pmod{p^n}
\]
for all sufficiently large primes $p$.
Equation \eqref{eq:rosen_kernel_equality} asserts that this congruence extends to a relation in $\widehat{\calA}$ without changing the coefficients of indices of weight less than $n$.
The homogeneous formulation in \cref{conj:Rosen_lifting} is equivalent to this one by multiplying or dividing by a suitable power of $\bp$ and adjusting the modulus accordingly.
\begin{proof}[Proof of \cref{thm:main} \ref{it:rosen_lifting}]
Assume \ref{it:dim_FMZV}.
We use the construction of $\calR_k^{(n)}$ also for $k=0$; the basis argument in \cref{prop:relations} applies without change.
For each positive integer $n$, there is a linear isomorphism $\alpha_n\colon\widehat{\calH}^1/\calI_n\xrightarrow{\sim}\calH_0^{(n)}$ with $\sfw(\bk)\bmod\calI_n\longmapsto t^{\lvert\bk\rvert}\sfw(\bk)$ for $\lvert\bk\rvert<n$.
Under this identification, the map induced by $Z_{\calA_n}$ is $Z_{\calA_n,0}$.
By \cref{prop:finite_injectivity}, $\psi_{\rKZ}^{\rcf,(n)}$ is injective on $\overline{\calZ^{\rcf}_{\calS_n,0}}$.
Thus \cref{lem:two-composite-rep,prop:relations} give $\ker Z_{\calA_n,0}=\ker Z_{\calS_n,0}^{\rcf}=\calR_0^{(n)}$.
Consequently,
\begin{equation}\label{eq:weighted_kernel}
\alpha_n\bigl(\ker Z_{\calA_n}/\calI_n\bigr)=\calR_0^{(n)}.
\end{equation}
We now lift the basis elements of $\calR_0^{(n)}$.
For $\bk\in\bbI\setminus\bbI^{\rZ}$, define
\begin{equation}\label{eq:rosen_generator}
\rho_{\bk}\coloneqq
\sfw(\bk)
-
\sum_{i=0}^{\infty}
\sum_{\bl\in\bbI^{\rZ}_{\lvert\bk\rvert+i}}
c_{\bk,\bl}^{(i)}\sfw(\bl)
\in\widehat{\calH}^1,
\end{equation}
where the coefficients are those of
\eqref{eq:infinite_formal_expansion}.
This series is well defined, since the inner sum is finite and all its words have length $\lvert\bk\rvert+i$.
For every positive integer $N$, the definition of $r_{j,\bk}^{(N)}$ gives
\[
\alpha_N(\rho_{\bk}\bmod\calI_N)=
\begin{cases}
r_{\lvert\bk\rvert,\bk}^{(N)}
&\text{if }\lvert\bk\rvert<N,\\
0&\text{if }\lvert\bk\rvert\geq N,
\end{cases}
\]
where the relation in the first case belongs to $\calR_0^{(N)}$.
Since $\calR_0^{(N)}\subset\ker Z_{\calA_N,0}$, we obtain $Z_{\calA_N}(\rho_{\bk})=0$ for every $N$.
It follows that $\rho_{\bk}\in\ker Z_{\widehat{\calA}}$.

Let
\[
r=\sum_{\bk\in\bbI}a_{\bk}\sfw(\bk)\in\ker Z_{\calA_n}.
\]
By \eqref{eq:weighted_kernel}, $\alpha_n(r\bmod\calI_n)$ belongs to $\calR_0^{(n)}$.
In its expansion in the basis given by \cref{prop:relations}, the coefficient of $r_{\lvert\bk\rvert,\bk}^{(n)}$ is $a_{\bk}$ for every
$\bk\in\bbI\setminus\bbI^{\rZ}$ with $\lvert\bk\rvert<n$.
Indeed, this basis element has coefficient $1$ at $t^{\lvert\bk\rvert}\sfw(\bk)$, and all its remaining terms have Zhao indices.
Hence
\[
\alpha_n(r\bmod\calI_n)=\sum_{\bk\in\bbI\setminus\bbI^{\rZ}, \ \lvert\bk\rvert<n}a_{\bk}r_{\lvert\bk\rvert,\bk}^{(n)}.
\]
Define the finite sum
\[
\widehat{r}\coloneqq\sum_{\bk\in\bbI\setminus\bbI^{\rZ}, \ \lvert\bk\rvert<n}a_{\bk}\rho_{\bk}.
\]
Each $\rho_{\bk}$ belongs to $\ker Z_{\widehat{\calA}}$, so the same is true of $\widehat{r}$.
Moreover,
\[
\alpha_n(\widehat{r}\bmod\calI_n)=\alpha_n(r\bmod\calI_n).
\]
Since $\alpha_n$ is injective, we have $r-\widehat{r}\in\calI_n$.
This proves
\[
\ker Z_{\calA_n}\subset\ker Z_{\widehat{\calA}}+\calI_n.
\]
The reverse inclusion follows immediately from reduction to $\calA_n$, proving \eqref{eq:rosen_kernel_equality}.

Finally, to compute $\widehat{r}\bmod\calI_N$ for any prescribed positive integer $N$, it suffices to compute the coefficients
$c_{\bk,\bl}^{(i)}$ occurring in \eqref{eq:rosen_generator} with
$\lvert\bk\rvert<n$ and $\lvert\bk\rvert+i<N$.
There are only finitely many such coefficients, and they are computed by the deterministic algorithm described in \cref{subsec:construction_of_rel}.
\end{proof}
%--------------------------------------------------------------------------
\section{Comparison of $\mathsf{FMZ}_{k,n}$ and $\calZ_{\calA_n,k}$}\label{sec:comparison}
%--------------------------------------------------------------------------
For a positive integer $n$ and a nonnegative integer $k$, set
\[
\mathsf{FMZ}_{k,n}\coloneqq\Span_{\bbQ}\{\zeta^{}_{\calA_n}(\bk):\bk\in\bbI_k\}.
\]
Then we have
\[
\calZ_{\calA_n,k}=\sum_{j=0}^{n-1}\bp^j\mathsf{FMZ}_{k+j,n}.
\]
In particular, $\mathsf{FMZ}_{k,n}\subset\calZ_{\calA_n,k}$ always holds, and it is natural to ask how these spaces differ.
The conjectural dimensions listed in \cite[Table~6]{Zhao} are interesting in this regard: they suggest that these spaces coincide for $n=2,3$, at least for $k\leq12$.
In this section, we prove that this equality indeed holds for $n=2,3$ and every nonnegative integer $k$.

By abuse of notation, we denote the restriction of $Z_{\calA_n,k}$ to $\calH^1_k$ by $\zeta^{}_{\calA_n}$.
Let $\be_j$ denote the vector whose $j$th component is $1$ and whose other components are $0$.
\begin{lemma}
For $\bk=(k_1,\dots, k_r)\in\bbI_k$ and $1\leq j\leq r$, we have
\begin{equation}\label{eq:DSR_2}
\bp k_j\zeta^{}_{\calA_2}(\bk+\be_j)\equiv 0 \pmod{\mathsf{FMZ}_{k,2}}
\end{equation}
and
\begin{equation}\label{eq:DSR_3}
\bp k_j\zeta^{}_{\calA_3}(\bk+\be_j)+\bp^2\binom{k_j+1}{2}\zeta^{}_{\calA_3}(\bk+2\be_j)+\bp^2\sum_{i=j+1}^rk_ik_j\zeta^{}_{\calA_3}(\bk+\be_i+\be_j)\equiv 0\pmod{\mathsf{FMZ}_{k,3}}.
\end{equation}
\end{lemma}
\begin{proof}
By the $\bp$-adic double shuffle relation first proved by Jarossay \cite{JarossayAdjoint} (cf.~\cite[Theorem~2.4]{OnoSakuradaSeki2021}), we have, for every positive integer $n$ and $0\leq i\leq r$,
\[
\sum_{\substack{\bl=(l_1,\dots, l_r)\in\bbZ^r_{\geq 0} \\ l_1=\cdots =l_i=0, \ \lvert\bl\rvert<n}}b(\bk,\bl)\bp^{\lvert\bl\rvert}\zeta^{}_{\calA_n}(\bk+\bl)=(-1)^{\lvert\bk^{[i]}\rvert}\zeta^{}_{\calA_n}(\sfw(\bk_{[i]})\sh \sfw(\overleftarrow{\bk^{[i]}}))\equiv 0 \pmod{\mathsf{FMZ}_{k,n}}.
\]
It suffices to subtract this identity with $i=j$ from that with $i=j-1$ and consider the cases $n=2,3$.
\end{proof}
%--------------------------------------------------------------
\subsection{The case $n=2$}
%--------------------------------------------------------------
\begin{theorem}
For every nonnegative integer $k$, we have
\[
\mathsf{FMZ}_{k,2}=\calZ_{\calA_2,k}.
\]
\end{theorem}
\begin{proof}
It suffices to show that $\bp\mathsf{FMZ}_{k+1,2}\subset\mathsf{FMZ}_{k,2}$
for every $k\geq0$.
Let $\bh=(h_1,\dots,h_r)\in\bbI_{k+1}$ be arbitrary.
We must show that $\bp\zeta^{}_{\calA_2}(\bh)\in\mathsf{FMZ}_{k,2}$.

Suppose that $h_j\geq2$ for some $j$, and choose $\bk\in\bbI_k$ such that $\bh=\bk+\be_j$.
By \eqref{eq:DSR_2}, we have $\bp(h_j-1)\zeta^{}_{\calA_2}(\bh)\in\mathsf{FMZ}_{k,2}$, which gives the desired inclusion since $h_j-1\neq0$.

In the remaining case, the well-known identity $\zeta^{}_{\calA}(\{1\}^{k+1})=0$ implies $\bp\zeta^{}_{\calA_2}(\{1\}^{k+1})=0$, completing the proof.
\end{proof}
%--------------------------------------------------------------
\subsection{The case $n=3$}
%--------------------------------------------------------------
\begin{theorem}
For every nonnegative integer $k$, we have
\[
\mathsf{FMZ}_{k,3}=\calZ_{\calA_3,k}.
\]
\end{theorem}
\begin{proof}
For $k=0$, we have
\[
\calZ_{\calA_3,0}=\bbQ+\bp\mathsf{FMZ}_{1,3}+\bp^2\mathsf{FMZ}_{2,3}.
\]
Wolstenholme's theorem gives $\zeta^{}_{\calA_2}(1)=0$, and hence $\bp\mathsf{FMZ}_{1,3}=0$.
Since $\zeta^{}_{\calA}(2)=\zeta^{}_{\calA}(1,1)=0$, we also have $\bp^2\mathsf{FMZ}_{2,3}=0$.
Thus, we may assume $k\geq1$ in what follows.

We first show that $\bp\mathsf{FMZ}_{k+1,3}\subset\bp^2\mathsf{FMZ}_{k+2,3}+\mathsf{FMZ}_{k,3}$.
Let $\bh=(h_1,\dots,h_r)\in\bbI_{k+1}$ and suppose that $h_j\geq2$ for some $j$.
Choose $\bk\in\bbI_k$ such that $\bh=\bk+\be_j$.
Then \eqref{eq:DSR_3} gives
\begin{equation}\label{eq:DSR_used1}
\bp\zeta^{}_{\calA_3}(\bh)\equiv-\bp^2\frac{h_j}{2}\zeta^{}_{\calA_3}(\bh+\be_j)-\bp^2\sum_{i=j+1}^r h_i\zeta^{}_{\calA_3}(\bh+\be_i)\pmod{\mathsf{FMZ}_{k,3}}.
\end{equation}
In particular, if $h_r\geq2$, then
\begin{equation}\label{eq:DSR_used2}
\bp\zeta^{}_{\calA_3}(\bh)\equiv-\bp^2\frac{h_r}{2}\zeta^{}_{\calA_3}(\sfw(\bh)e_0)\pmod{\mathsf{FMZ}_{k,3}}.
\end{equation}
Define an automorphism $\psi$ of $\calH$ by $\psi(e_0)=e_0+e_1$ and $\psi(e_1)=-e_1$.
Then Rosen's $\bp$-adic duality \cite[Theorem~4.5]{Rosen2015} (cf.\ \cite[Theorem~5.1]{MaesakaSekiWatanabe}) gives
\begin{equation}\label{eq:Rosen_duality}
\bp\zeta^{}_{\calA_3}(we_1)-\bp^2\zeta^{}_{\calA_3}((w\ast e_1)e_1)=\zeta^{}_{\calA_3}(\psi(w))-\zeta^{}_{\calA_3}(w)\equiv0\pmod{\mathsf{FMZ}_{k,3}}
\end{equation}
for every nonempty word $w=\sfw(\bk)\in e_1\calH$ with $\bk\in\bbI_k$, where $\ast$ denotes the stuffle product.
Consequently, if $\sfw(\bh)=\sfw(\bk)e_1$ for $\bh\in\bbI_{k+1}$, then $\bp\zeta^{}_{\calA_3}(\bh)\equiv\bp^2\zeta^{}_{\calA_3}((w\ast e_1)e_1)\pmod{\mathsf{FMZ}_{k,3}}$.
Since $(we_1)\ast e_1=(w\ast e_1)e_1+we_1^2+we_1e_0$ and $\zeta^{}_{\calA}(1)=0$, we have $\zeta^{}_{\calA}((w\ast e_1)e_1)=-\zeta^{}_{\calA}(\sfw(\bh)e_1)-\zeta^{}_{\calA}(\sfw(\bh)e_0)$.
Hence we obtain
\begin{equation}\label{eq:Rosen_duality_used}\bp\zeta^{}_{\calA_3}(\bh)\equiv-\bp^2\zeta^{}_{\calA_3}(\sfw(\bh)e_1)-\bp^2\zeta^{}_{\calA_3}(\sfw(\bh)e_0)\pmod{\mathsf{FMZ}_{k,3}}.
\end{equation}
Equations \eqref{eq:DSR_used2} and \eqref{eq:Rosen_duality_used} therefore establish $\bp\mathsf{FMZ}_{k+1,3}\subset\bp^2\mathsf{FMZ}_{k+2,3}+\mathsf{FMZ}_{k,3}$.
Thus, it remains to show that $\bp^2\mathsf{FMZ}_{k+2,3}\subset\mathsf{FMZ}_{k,3}$.

We now show that, for every $\bl\in\bbI_{k-1}$,
\begin{equation}\label{eq:21-111}
2\bp^2\zeta^{}_{\calA_3}(\sfw(\bl)e_1e_0e_1)+\bp^2\zeta^{}_{\calA_3}(\sfw(\bl)e_1^3)\equiv0\pmod{\mathsf{FMZ}_{k,3}}.
\end{equation}
The equality in \eqref{eq:Rosen_duality} holds for words of arbitrary weight and extends by linearity to $e_1\calH$.
Multiplying it by $\bp$ and $\bp^2$, respectively, we obtain
\begin{equation}\label{eq:Rosen_duality_p}
\bp\zeta^{}_{\calA_3}(\psi(w))-\bp\zeta^{}_{\calA_3}(w)=\bp^2\zeta^{}_{\calA_3}(we_1)
\end{equation}
and
\begin{equation}\label{eq:Rosen_duality_p^2}
\bp^2\zeta^{}_{\calA_3}(\psi(w))=\bp^2\zeta^{}_{\calA_3}(w).
\end{equation}
Substituting $w=\sfw(\bl)e_1e_0$ into
\eqref{eq:Rosen_duality_p} gives $\bp\zeta^{}_{\calA_3}(\psi(\sfw(\bl)e_1e_0))-\bp\zeta^{}_{\calA_3}(\sfw(\bl)e_1e_0)=\bp^2\zeta^{}_{\calA_3}(\sfw(\bl)e_1e_0e_1)$.
By \eqref{eq:DSR_used2}, we have
\begin{equation}\label{eq:p_equiv_-p^2}
\bp\zeta^{}_{\calA_3}(\sfw(\bl)e_1e_0)\equiv-\bp^2\zeta^{}_{\calA_3}(\sfw(\bl)e_1e_0^2)\pmod{\mathsf{FMZ}_{k,3}}.
\end{equation}
Since $\psi(\sfw(\bl)e_1e_0)=-\psi(\sfw(\bl))e_1e_0-\psi(\sfw(\bl))e_1^2$, applying \eqref{eq:DSR_used2} and \eqref{eq:Rosen_duality_used} by linearity yields
\[
\bp\zeta^{}_{\calA_3}(\psi(\sfw(\bl)e_1e_0))\equiv\bp^2\zeta^{}_{\calA_3}\bigl(\psi(\sfw(\bl))(e_1e_0^2+e_1^3+e_1^2e_0)\bigr)\pmod{\mathsf{FMZ}_{k,3}}.
\]
By \eqref{eq:Rosen_duality_p^2} and $\psi^2=\mathrm{id}$,
the right-hand side equals $\bp^2\zeta^{}_{\calA_3}\bigl(\sfw(\bl)\psi(e_1e_0^2+e_1^3+e_1^2e_0)\bigr)$.
Since $\psi(e_1e_0^2+e_1^3+e_1^2e_0)=-e_1e_0^2-e_1e_0e_1-e_1^3$,
we can eliminate
$\bp\zeta^{}_{\calA_3}(\psi(\sfw(\bl)e_1e_0))$
and use \eqref{eq:p_equiv_-p^2} to obtain \eqref{eq:21-111}.

Suppose that $\bh\in\bbI_{k+1}$ satisfies $\sfw(\bh)=\sfw(\bv)e_1e_0^{b-1}e_1$ with $b\geq2$.
Eliminating $\bp^2\zeta^{}_{\calA_3}(\sfw(\bh)e_0)$ between the congruences obtained from \eqref{eq:DSR_used1} and \eqref{eq:Rosen_duality_used} gives
\begin{equation}\label{eq:1,1_+1,1}
\bp^2\zeta^{}_{\calA_3}(\sfw(\bv)e_1e_0^{b-1}e_1^2)\equiv\bp^2\frac{b}{2}\zeta^{}_{\calA_3}(\sfw(\bv)e_1e_0^be_1)\pmod{\mathsf{FMZ}_{k,3}}.
\end{equation}
Next, suppose that $\bh\in\bbI_{k+1}$ satisfies $\sfw(\bh)=\sfw(\bv)e_1e_0^{b-1}e_1^2$ with $b\geq2$.
Using \eqref{eq:DSR_used1} and \eqref{eq:Rosen_duality_used}, we obtain
\[
\bp^2\zeta^{}_{\calA_3}(\sfw(\bv)e_1e_0^{b-1}e_1^3)\equiv\bp^2\frac{b}{2}\zeta^{}_{\calA_3}(\sfw(\bv)e_1e_0^be_1^2)+\bp^2\zeta^{}_{\calA_3}(\sfw(\bv)e_1e_0^{b-1}e_1e_0e_1)\pmod{\mathsf{FMZ}_{k,3}}.
\]
Applying \eqref{eq:21-111} to the left-hand side and \eqref{eq:1,1_+1,1} to the first term on the right-hand side gives
\begin{align*}
&-2\bp^2\zeta^{}_{\calA_3}(\sfw(\bv)e_1e_0^{b-1}e_1e_0e_1)\\
&\equiv\frac{b(b+1)}{4}\bp^2\zeta^{}_{\calA_3}(\sfw(\bv)e_1e_0^{b+1}e_1)+\bp^2\zeta^{}_{\calA_3}(\sfw(\bv)e_1e_0^{b-1}e_1e_0e_1)\pmod{\mathsf{FMZ}_{k,3}}.
\end{align*}
Therefore, we have
\begin{equation}\label{eq:b(b+1)/12}
\bp^2\zeta^{}_{\calA_3}(\sfw(\bv)e_1e_0^{b-1}e_1e_0e_1)\equiv-\frac{b(b+1)}{12}\bp^2\zeta^{}_{\calA_3}(\sfw(\bv)e_1e_0^{b+1}e_1)\pmod{\mathsf{FMZ}_{k,3}}.
\end{equation}
Let $b\geq2$ and $c\geq3$, and apply \eqref{eq:DSR_used1} to $\bh\in\bbI_{k+1}$ satisfying $\sfw(\bh)=\sfw(\bv)e_1e_0^{b-1}e_1e_0^{c-2}e_1$.
Comparing the two congruences obtained at the positions with entries $b$ and $c-1$, we obtain
\[
(c-1)\bp^2\zeta^{}_{\calA_3}(\sfw(\bv)e_1e_0^{b-1}e_1e_0^{c-1}e_1)\equiv-b\bp^2\zeta^{}_{\calA_3}(\sfw(\bv)e_1e_0^be_1e_0^{c-2}e_1)\pmod{\mathsf{FMZ}_{k,3}}.
\]
Repeating this step $c-2$ times and finally applying \eqref{eq:b(b+1)/12}, we obtain
\begin{equation}\label{eq:q_b,c}
\bp^2\zeta^{}_{\calA_3}(\sfw(\bv)e_1e_0^{b-1}e_1e_0^{c-1}e_1)\equiv q_{b,c}\bp^2\zeta^{}_{\calA_3}(\sfw(\bv)e_1e_0^{b+c-1}e_1)\pmod{\mathsf{FMZ}_{k,3}}
\end{equation}
for some rational number $q_{b,c}$.
This also holds for $c=2$ by \eqref{eq:b(b+1)/12}, with $q_{b,2}=-b(b+1)/12$.
All factors arising in the iteration are nonzero, so $q_{b,c}\neq0$ for all $b,c\geq2$.

By \cref{thm:previous} \eqref{it:FMZV_basis}, to prove $\bp^2\mathsf{FMZ}_{k+2,3}\subset\mathsf{FMZ}_{k,3}$, it suffices to show that $\bp^2\zeta^{}_{\calA_3}(\bl)\in\mathsf{FMZ}_{k,3}$ for every $\bl\in\bbI^{\rZ}_{k+2}$.
Repeatedly applying \eqref{eq:q_b,c} to reduce the depth, we find that $\bp^2\zeta^{}_{\calA_3}(\bl)$ is congruent to a nonzero rational multiple of $\bp^2\zeta^{}_{\calA_3}(k+1,1)$ modulo $\mathsf{FMZ}_{k,3}$.

If $k$ is even, then $\zeta^{}_{\calA}(k+1,1)=0$, and hence $\bp^2\zeta^{}_{\calA_3}(k+1,1)=0$.
If $k$ is odd, we have
\[
\bp^2\zeta^{}_{\calA_3}(k+1,1)=\frac{2(k+2)}{k(k+1)}\zeta^{}_{\calA_3}(k)\in\mathsf{FMZ}_{k,3}.
\]
For these special values, see \cite[Section~3]{OnoSakuradaSeki2021}.
This completes the proof.
\end{proof}
Establishing the four assumptions \ref{it:dim_MZV}, \ref{it:dim_FMZV}, \ref{it:direct_sum_MZV}, and \ref{it:direct_sum_FMZV} remains a major goal for future research. For ordinary multiple zeta values, Fauzan's proof of the irrationality of $\zeta(5)$ in 2026 \cite{Fauzan} marked a significant advance, but we are still far from proving either \ref{it:dim_MZV} or \ref{it:direct_sum_MZV}. For finite multiple zeta values, even the existence of a nonempty index $\bk$ for which $\zeta^{}_{\calA}(\bk)\ne 0$ remains an open problem. As advances in generative AI transform mathematics, when will humanity finally succeed in proving these conjectures?
%--------------------------------------------------------------
\subsection*{Use of AI}
%--------------------------------------------------------------
ChatGPT (OpenAI, GPT-6 Astra Pro) was used for mathematical discussions and English-language editing.
The final manuscript and all mathematical arguments were written and independently verified by the author, who takes full responsibility for the content.
%--------------------------------------------------------------

\end{document}